\documentclass[12pt]{article}

\usepackage{amsxtra,amssymb,amsthm,amsmath,latexsym}
\usepackage{graphicx}
\usepackage{epsfig}
\usepackage{epstopdf}
\usepackage{booktabs}
\usepackage{multirow}
\usepackage{afterpage}
\usepackage{setspace}
\usepackage{microtype}
\usepackage{calc,fourier}
\usepackage[sort&compress]{natbib}
\usepackage[T1]{fontenc}
\usepackage[hidelinks]{hyperref}
\usepackage{float}
\usepackage[small]{caption}
\newtheorem{thm}{Theorem}[]

\numberwithin{equation}{section}

\newcommand{\bee}{\begin{equation*}}
\newcommand{\eee}{\end{equation*}}
\newcommand{\be}{\begin{equation}}
\newcommand{\ee}{\end{equation}}
\newcommand{\ba}{\begin{align}}
\newcommand{\ea}{\end{align}}

\newcommand{\RRR}{\mathbb{R}^3}

\title{Numerical Solution of Many-body Wave Scattering Problem and Creating Materials with A Desired Refraction Coefficient}
\author{N. T. Tran\footnote{Mailing address:  Mathematics Department, 138 Cardwell Hall, Manhattan, KS 66506} \\
\small Department of Mathematics\\[-0.8ex]
\small Kansas State University, Manhattan, KS 66506-2602, USA\\
\small \texttt{nhantran@math.ksu.edu}}

\date{}
\begin{document}
\maketitle

\begin{abstract}
Scalar wave scattering by many small particles with impedance boundary condition and creating material with a desired refraction coefficient are studied. The acoustic wave scattering problem is solved asymptotically and numerically under the assumptions $ka  \ll  1, \zeta_m = \frac{h(x_m)}{a^\kappa}, d = O(a^{\frac{2-\kappa}{3}}), M = O(\frac{1}{a^{2-\kappa}}), \kappa \in [0,1)$, where $k = 2\pi/\lambda$ is the wave number, $\lambda$ is the wave length, $a$ is the radius of the particles, $d$ is the distance between neighboring particles, $M$ is the total number of the particles embedded in a bounded domain $\Omega \subset \RRR$, $\zeta_m$ is the boundary impedance of the m\textsuperscript{th} particle $D_m$, $h \in C(D)$, $D := \bigcup_{m=1}^M D_m$, is a given arbitrary function which satisfies Im$h \le 0$, $x_m \in \Omega$ is the position of the m\textsuperscript{th} particle, and $1 \leq m \leq M$. Numerical results are presented for which the number of particles equals $10^4, 10^5$, and $10^6$.
\end{abstract}



\section{Introduction} \label{sec0}
Recent research in materials science shows the existence of materials with negative refraction coefficient, called metamaterials; see \citet{EB2005}. Creating such materials is of practical interest since metamaterials have many applications but are not available in nature; e.g. see \citet{Hansen2008}. By arranging their structure, one can create new materials with a desired refraction coefficient. In \citet{R13,R509,R536,R563,R580,R586,R591,R595,R632,R635}, A. G. Ramm has developed the theory of wave scattering by many small bodies for acoustic and electromagnetic (EM) waves that can be used for creating materials with a desired refraction coefficient.

In \citet{R595,R635}, he derived analytic formulas for the solution of wave scattering by many small bodies (ie), the algebraic system (ori), and the reduced order system (red) for asymptotically solving wave scattering by many small bodies of arbitrary shapes, and developed an approach for creating materials with a desired refraction coefficient. The small bodies can be particles whose physical properties are described by their boundary impedance. This theory can be used in many practical problems. In \citet{R595,R635}, an asymptotic solution of the many-body acoustic wave scattering problem was developed under the assumptions $ka  \ll  1, \zeta_m = \frac{h(x_m)}{a^\kappa}, d = O(a^{\frac{2-\kappa}{3}}), M = O(\frac{1}{a^{2-\kappa}})$, and $\kappa \in [0,1)$, where $k = 2\pi/\lambda$ is the wave number, $a$ is the radius of the particles, defined as $a:=\frac{1}{2}\max_{1 \le m \le M}\text{diam}D_m$ and $D_m$ is the m\textsuperscript{th} particle, $d$ is the distance between neighboring particles, $M$ is the total number of the particles embedded in a bounded domain $\Omega \subset \RRR$, $\zeta_m$ is the boundary impedance of the m\textsuperscript{th} particle, $h \in C(D)$, $D := \bigcup_{m=1}^M D_m$, is a given arbitrary continuous function, Im $h \leq 0$, $x_m \in D_m$ is an arbitrary point in the m\textsuperscript{th} particle, and $1 \leq m \leq M$.

This paper will focus on solving the systems (ori), (red), and integral equation (ie) for wave scattering problem by many small impedance particles with complex refraction coefficients. The goal is to check the numerical accuracy of the solutions to (ori), (red), and (ie) when the number of particles is large, up to order $10^6$. There was no results on solving wave scattering problem for so many particles as in this paper. Furthermore, these results are used for creating materials with a desired refraction coefficient, as was proved in \cite{R635}. In this paper, the theory from \cite{R635} is illustrated by numerical examples.

\section{Wave scattering by one small impedance particle} \label{sec1}
Let us formulate the wave scattering problem with one body. Let $D$ be a bounded domain of one small particle in $\RRR$ , $D'$ be the exterior domain of $D$, and $S$ be the boundary of $D$. Let $\alpha \in S^2$ denote the direction of the incident plane wave, $|\alpha|=1$, and $S^2$ denote a unit sphere. Finally, let $u_0$ be the incident field that satisfies Helmholtz equation in $\RRR$, $v$ be the scattered field which satisfies the radiation condition, and $a$ be the radius of the particle. Then the scattering problem consists of solving the following system:
\begin{align}
    &(\nabla^2+k^2)u(x)=0 \quad\text{ in } D',  \quad k=\text{const}>0, \label{eq1.1}\\
    &u_N=\zeta u \quad \text{ on }S, \quad \text{Im }\zeta \leq 0, \label{eq1.2}\\
    &u(x) = u_0(x)+v(x),  \label{eq1.3} \\
    &u_0(x) = e^{ik\alpha \cdot x}, \quad ka  \ll  1, \quad\text{and $v$ satisfies the radiation condition:} \label{eq1.4} \\
    &v_r-ikv=o(1/r), \quad r:=|x| \to \infty, \label{eq1.5}
\end{align}
where $k$ is a wave number, $\zeta$ is the boundary impedance of the surface $S$, and $N$ is the outer unit normal vector to $S$.
If Im$\zeta \le 0$, it was proved in \citet{R632} that the system \eqref{eq1.1}-\eqref{eq1.5} has a unique solution of the form
\be
    u(x)=u_0(x)+\int_S g(x,t)\sigma(t)dt,
\ee
where $g(x,t):=\frac{e^{ik|x-t|}}{4\pi |x-t|}$ and $\sigma(t)$ is some continuous function which is uniquely defined by the boundary condition.

The asymptotic solution is then given by the formula
\be
    u \sim u_0+ g(x,x_1)Q, \quad \text{when } |x| \gg a, \quad a \to 0.
\ee
Here
\be
    Q:=\int_S \sigma(t)dt \simeq -\zeta|S|u_0(x_1), \quad a \to 0,
\ee
where $|S|$ is the surface area of the small particle $D$, and $x_1 \in D$. Instead of finding function $\sigma(t)$ to get the solution $u$, one can just find the number $Q$.

\section{Wave scattering by many small impedance particles} \label{sec2}
Consider a bounded domain $\Omega \subset \RRR$ that is filled with a material consisting of $M$ particles. This material has refraction coefficient $n_0(x)$. Let $D_m$ be the domain of one particle and $S_m$ be the boundary of $D_m$. Define $D := \bigcup_{m=1}^M D_m \subset \Omega$ and $D':=\RRR \setminus D$. The minimal distance between neighboring particles, $d$, is much greater than the maximal radius of a particle, $a=\frac{1}{2}\max_{1 \leq m \leq M} \text{diam}D_m$, and much less than $\lambda$, the wave length. Let $\zeta_m$ denote the boundary impedance of $S_m$, $\zeta_m=\frac{h(x_m)}{a^\kappa}$, where $h(x)$ is a continuous function in $D$ such that Im $h \leq 0$ in $D$ and $\kappa$ is a const in [0,1). The scattering problem is then formulated as follows:
\begin{align}
    &(\nabla^2+k^2 n_0^2(x))u=0 \quad\text{ in } D', \quad k=\text{const}>0, \label{eq2.1} \\
    &u_N=\zeta_m u \quad \text{ on }S_m, \quad \text{Im }\zeta_m \leq 0,  \quad 1 \leq m \leq M, \label{eq2.2} \\
    &u(x) = u_0(x)+v(x),  \label{eq2.3} \\
    &u_0(x) = e^{ik\alpha \cdot x}, \quad ka  \ll  1, \quad\text{and $v$ satisfies the radiation condition:} \label{eq2.4} \\
    &v_r-ikv=o(1/r), \quad r:=|x| \to \infty, \label{eq2.5}
\end{align}
where $k$ is a wave number and $n_0(x)=1$ in $\Omega'$ is the initial refraction coefficient such that Im $n_0^2(x) \geq 0$ in $\Omega$ and it is a Riemann-integrable function. It was proved in \citet{R536} that if Im $n_0^2(x) \geq 0$ and Im $h(x) \leq 0$, then the system \eqref{eq2.1}-\eqref{eq2.5} has a unique solution of the form
\be \label{eq3}
    u(x)=u_0(x)+\sum_{m=1}^M \int_{S_m} G(x,y)\sigma_m(y)dy,
\ee
where $G(x,y)$ is a Green function of the Helmholtz equation \eqref{eq2.1}, $G$ satisfies $[\nabla^2 + k^2 n_0^2(x)] G = -\delta(x-y)$ in $\RRR$ and the radiation condition, and $\sigma_m(y)$ are some continuous functions which are uniquely defined by the boundary condition.

Let us assume for simplicity that $x_j$ is the center of $D_j$, a ball of radius $a$. Then we define the effective field acting on the j\textsuperscript{th} particle as
\be \label{eq4}
    u_e(x_j):=u(x)-\int_{S_j} G(x_j,y)\sigma_j(y)dy,
\ee
or equivalently
\be \label{eq5}
    u_e(x_j)=u_0(x_j)+\sum_{m=1 ,m \neq j}^M \int_{S_m} G(x_j,t)\sigma_m(t)dt.
\ee
Let us derive the approximation formula for this effective field. From \eqref{eq3}, one gets
\be \label{eq6}
    u(x)=u_0(x) + \sum_{m=1}^M G(x,x_m)Q_m + \sum_{m=1}^M \int_{S_m} [G(x,y)-G(x,x_m)] \sigma_m(y)dy.
\ee
Here
\be
    Q_m:=\int_{S_m} \sigma_m(y)dy.
\ee
Instead of finding functions $\sigma_m(y)$ to get the solution $u$, one can just find numbers $Q_m$.

One can rewrite \eqref{eq6} as follows (see \citet{R509}, \citet{R536}):
\be \label{eq7}
    u(x)=u_0(x)+\sum_{m=1}^M G(x,x_m)Q_m + o(1),
\ee
as $a \to 0$ and  $|x-x_m| \geq a$. When $a \to 0$, one can compute $Q_m$ asymptotically and get
\be \label{eq8}
    Q_m \simeq -c a^{2-\kappa} h(x_m)u_e(x_m),
\ee
where $c$ is a constant depending on the shape of a particle, $|S|=ca^2$, where $|S|$ is the surface area of $S$. If $S$ is a sphere, then $c=4\pi$.
Thus, one can rewrite \eqref{eq5} as
\be \label{eq10}
    u_e(x_j) \simeq u_0(x_j)-4\pi\sum_{m=1, m \neq j}^M G(x_j,x_m)h(x_m)u_e(x_m)a^{2-\kappa},
\ee
as $a \to 0$ and $1 \leq j \leq M$.
Denote $u_j:=u_e(x_j),  u_{0j}:=u_0(x_j), G_{jm}:=G(x_j,x_m)$, and $h_m:=h(x_m)$. In \eqref{eq10}, the numbers $u_m$, $1 \le m \le M$, are unknowns. It was proved in \citet{R595,R635} that under the assumptions
\be \label{eq11}
    d=O\left(a^{\frac{2-\kappa}{3}}\right), \quad\text{and } M=O\left(\frac{1}{a^{2-\kappa}}\right), \quad\text{for } \kappa \in [0,1),
\ee
$u_j$, where $1 \leq j \leq M,$ can be found by solving the linear algebraic system (LAS)
\be \label{eq12}
    u_j = u_{0j}-4\pi\sum_{m=1, m \neq j}^M G_{jm} h_m a^{2-\kappa} u_m, \quad\text{as } a \to 0, \quad 1 \leq j \leq M.
\ee
We call this LAS the original system (ori).

Let $\Delta$ be a subdomain in $\Omega$ and $\mathcal{N}(\Delta)$ be the number of embedded particles in $\Delta$. We assume that
\be \label{eq13}
    \mathcal{N}(\Delta)=\frac{1}{a^{2-\kappa}} \int_{\Delta} N(x)dx[1+o(1)], \quad\text{as } a \to 0,
\ee
where $N(x) \geq 0$ is a given continuous function in $\Omega$, $N(x)$ and $\kappa$ can be chosen as desired.

Let $\Omega$ be partitioned into $P$ non-intersecting sub cubes $\Delta_p$'s of size $b$ such that $b \gg d \gg a$, where $b=b(a), d=d(a),$ and $\lim_{a\to 0} \frac{d(a)}{b(a)}=0$. One can then derive, see \citet{R595,R635}, from \eqref{eq12} and \eqref{eq13} that
\be \label{eq14}
    u_q = u_{0q}-4\pi\sum_{p=1, p \neq q}^P G_{qp} h_p N_p u_p |\Delta_p|, \quad\text{for } 1 \leq q \leq P,
\ee
where $|\Delta_p|$ is the volume of $\Delta_p$, $N_p:=N(x_p)$, and $x_p$ is a point in $\Delta_p$, for example, the center of $\Delta_p$. This linear system is much easier to solve than (ori) since $P  \ll  M$. We will call the LAS \eqref{eq14} the reduced ordered system (red).

If assumption \eqref{eq13} holds, the limiting integral equation obtained from \eqref{eq14} as $a \to 0$ is
\be \label{eq15}
    u(x)=u_0(x)-4\pi\int_D G(x,y)h(y)N(y)u(y)dy, \quad\text{for } x \in \mathbb{R}^3,
\ee
or equivalently
\be \label{eq16}
    u(x)=u_0(x)-\int_D G(x,y)p(y)u(y)dy, \quad\text{for } x \in \mathbb{R}^3,
\ee
where $p(x):=4\pi h(x)N(x)$.
This integral equation yields the limiting field in the medium created by embedding many small particles with distribution \eqref{eq13}; see \citet{R536} and \citet{R595,R635}. Any function $p(x)$ can be created by choosing functions $h(x)$ and $N(x)$ properly; see Section \ref{sec3}. We will call equation \eqref{eq16} the integral equation (ie).

The following result was proved in \citet{R595,R635}.
\begin{thm}
If assumptions \eqref{eq11} and \eqref{eq13} hold, then there exists the limit
\be \label{eq17}
    \lim_{a \to 0} ||u_e(x)-u(x)||_{C(\mathbb{R}^3)}=0,
\ee
where $u(x)$ is the unique solution to  (ie).
\end{thm}

\section{A recipe for creating materials with a desired refraction coefficient} \label{sec3}
We want to create from the material with initial refraction coefficient $n_0(x)$ a new material with a desired refraction coefficient $n(x)$. We describe the recipe, proposed in \citet{R595,R635}, to accomplish this. This recipe has three steps.

{\bf\em Step 1:} Calculate p(x) using the following formula whose derivation can be found in \citet{R595,R635}
\be \label{eq18}
    p(x)=k^2[n_0^2(x)-n^2(x)].
\ee

{\bf\em Step 2:} Choose an arbitrary $N(x)>0$ and use the relation $p(x)=4\pi h(x)N(x)$ to calculate $h(x):=h_1(x)+i h_2(x)$ as follows
\be \label{eq19}
    h_1(x)=\frac{p_1(x)}{4\pi N(x)}, \quad h_2(x)=\frac{p_2(x)}{4\pi N(x)},
\ee
where $p_1(x)=$ Re $p(x)$ and $p_2(x)=$ Im $p(x)$. Note that Im $h(x) \leq 0$ holds if Im $p(x) \leq 0$.

{\bf\em Step 3:} Embed $M$ small particles of radius $a$ with boundary impedance $\zeta_m = \frac{h(x_m)}{a^\kappa}$, where $1 \leq m \leq M$ and $M=\frac{1}{a^{2-\kappa}}\int_\Omega N(x) dx [1+o(1)]$, into the domain $\Omega$ at the approximately prescribed positions according to formula \eqref{eq13}.

The resulting materials, obtained by embedding many small particles into $\Omega$ using this recipe, will have the desired refraction coefficient $n(x)$ with an error that tends to zero as $a \to 0$, as proved in \citet{R536}.

\section{Numerical results} \label{sec4}
In this section, we present some numerical results of solving the wave scattering problem by many small particles, in particular, solving  (ori), (red), and (ie). For solving linear algebraic systems, we used PETSC libraries developed at Argonne National Lab to do the computation in parallel; see \citet{Petsc}. GMRES iterative method, see \citet{Gmres}, is used to find the solutions to (ori) and (red) with relative error equal to $10^{-3}$. For solving (ie), we used the collocation method from \citet{R563}, dividing the domain into many sub cubes, taking the collocation points as the centers of these cubes, and then approximating the integral equation by the corresponding Riemann sum. After that, we used GMRES iterative method to find an approximation of the solution to (ie) with relative error equal to $10^{-3}$. Since the number of unknowns in (ori), (red), and (ie) are different, we used an interpolation procedure to compare their solutions. For example, let the domain $\Omega$ be a unit cube that contains $M$ particles. We partitioned $\Omega$ into $P$ small sub cubes to solve (red). In this case, (ori) has $M$ unknowns, say $x_i$, $1 \leq i \leq M$, and (red) has $P$ unknowns, say $y_q$ for $1 \leq q \leq P$. Let us assume that $M>P$. To find the difference between solutions to (ori) and (red), we find all the particles $x_i$ that lie in a sub cube $\Delta_q$ corresponding to $y_q$, and then find the solution differences $|x_i-y_q|$ for these particles. After that, we compute the following
\be \label{eq20}
    \sup_{y_q} \frac{1}{\mathcal{N}(\Delta_q)}\sum_{x_i \in \Delta_q} |x_i-y_q|,
\ee
where $\mathcal{N}(\Delta_q)$ is the number of particles in the sub cube $\Delta_q$. This gives us the solution difference between (ori) and (red). The solution differences between (ori) and (ie), and (ie) and (red) are computed similarly.

The following numerical experiments are of practical interest and importance. One wants to find:\\
a) The solution differences between  (ori) and (red), (ie) and (ori), and (ie) and (red), denoted $e_1, e_2$, and $e_3$, respectively. \\
b) The maximal value of $a/d$ for which the solution differences are less than 3\% or  5\%, for example. \\
c) The values of $a/d$ for which the solution difference becomes larger than say 10\%, i.e. for which the asymptotic formula \eqref{eq7} is no longer applicable.

The error considered later is the solution difference $e=e_1+e_2+e_3$. One can find in \citet{R580} numerical results for $M \le 15^3$ particles. In this paper, we will do the experiment with a large number $M$ of particles, such as $M=10^4, 10^5$ or $10^6$. We assume that the domain $\Omega$ that contains all the particles, is a unit cube. The following values of physical parameters are used to conduct the experiment:

- Wave number, $k=0.182651 \text{ cm}^{-1}$;

- Direction of the incident plane wave, $\alpha=$   (1, 0, 0);

- The constant $\kappa=$ 	0.99;

- Volume of the domain $\Omega$ that contains all the particles, $|\Omega|=1 \text{ cm}^3$;

- Original refraction coefficient, $n_0=$ 1+0i;

- Desired refraction coefficient, $n=$  -1+0.001i;

- The function $N(x)=Ma^{2-\kappa}/|\Omega|$, i.e. particles are distributed uniformly in the unit cube;

- Number of small sub cubes after partitioning the domain $\Omega$ for solving (red), $P=$ 125;

- Number of collocation points for solving (ie), $C=$ 8000.

In this case, $P  < C <  M$. To do the solution comparisons, the interpolation procedure described above is used to obtain the solutions to (red) and (ie) at the points corresponding to the position of the particles.

\begin{figure}[htbp]
    \centering
    \includegraphics[scale=1.1]{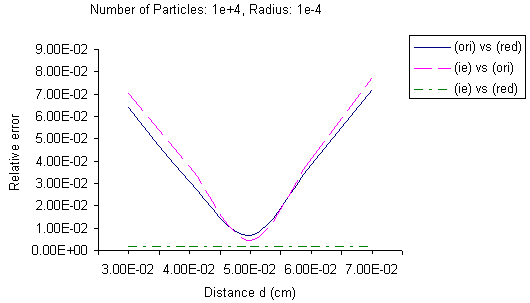}
    \caption{Solution comparison of (ori), (red), and (ie) with $M=10^4$, $a=10^{-4}$, and different $d$.}
    \label{fig1}
\end{figure}

\begin{figure}[htbp]
    \centering
    \includegraphics[scale=1.1]{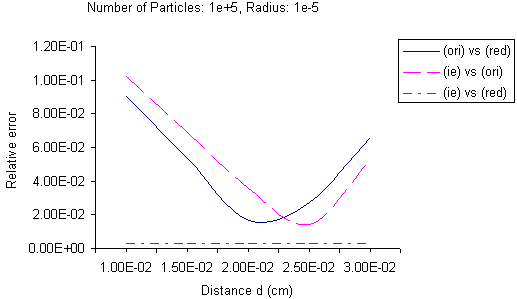}
    \caption{Solution comparison of (ori), (red), and (ie) with $M=10^5$, $a=10^{-5}$, and different $d$.}
    \label{fig2}
\end{figure}

\begin{figure}[htbp]
    \centering
    \includegraphics[scale=1.1]{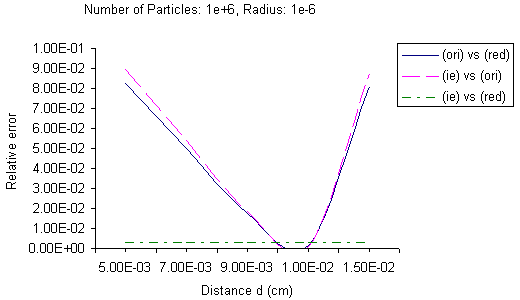}
    \caption{Solution comparison of (ori), (red), and (ie) with $M=10^6$, $a=10^{-6}$, and different $d$.}
    \label{fig3}
\end{figure}

In Figure \ref{fig1}, \ref{fig2}, and \ref{fig3}, the solid line shows the difference between solutions to  (ori) and (red), the dashed line shows the difference between solutions to  (ie) and (ori), and the dot-dashed line shows the difference between solutions to  (ie) and (red). Radius of particles and distance between neighboring particles are measured in centimeters. We will consider the error sum $e$, the sum of the three solution differences, to choose the best ratio $a/d$ for each $a$.
\begin{table}[H]
  \centering
  \tabcolsep=0.11cm
  \small
  \caption{Solution comparison of (ori), (red), and (ie) with $M=10^4$, $a=10^{-4}$, and different $d$.}
    \begin{tabular}{crrrrrr}
    \toprule
    \multicolumn{7}{c}{M=1.00E+4, a=1.00E-4}
      \\
    \midrule
    d &
      \multicolumn{1}{l}{2.00E-02} &
      \multicolumn{1}{l}{3.00E-02} &
      \multicolumn{1}{l}{4.00E-02} &
      \multicolumn{1}{l}{5.00E-02} &
      \multicolumn{1}{l}{6.00E-02} &
      \multicolumn{1}{l}{7.00E-02}
      \\
    a/d &
      \multicolumn{1}{l}{5.00E-03} &
      \multicolumn{1}{l}{3.33E-03} &
      \multicolumn{1}{l}{2.50E-03} &
      \multicolumn{1}{l}{2.00E-03} &
      \multicolumn{1}{l}{1.67E-03} &
      \multicolumn{1}{l}{1.43E-03}
      \\
    (ori) vs. (red) &
      \multicolumn{1}{l}{9.75E-02} &
      \multicolumn{1}{l}{6.41E-02} &
      \multicolumn{1}{l}{3.07E-02} &
      \multicolumn{1}{l}{6.52E-03} &
      \multicolumn{1}{l}{3.84E-02} &
      \multicolumn{1}{l}{7.21E-02}
      \\
    (ie) vs. (ori) &
      \multicolumn{1}{l}{1.05E-01} &
      \multicolumn{1}{l}{7.03E-02} &
      \multicolumn{1}{l}{3.74E-02} &
      \multicolumn{1}{l}{4.57E-03} &
      \multicolumn{1}{l}{4.11E-02} &
      \multicolumn{1}{l}{7.76E-02}
      \\
    (ie) vs. (red) &
      \multicolumn{1}{l}{1.83E-03} &
      \multicolumn{1}{l}{1.83E-03} &
      \multicolumn{1}{l}{1.83E-03} &
      \multicolumn{1}{l}{1.83E-03} &
      \multicolumn{1}{l}{1.83E-03} &
      \multicolumn{1}{l}{1.83E-03}
      \\
    Error sum $e$ &
      \multicolumn{1}{l}{2.04E-01} &
      \multicolumn{1}{l}{1.36E-01} &
      \multicolumn{1}{l}{7.00E-02} &
      \multicolumn{1}{l}{1.29E-02} &
      \multicolumn{1}{l}{8.13E-02} &
      \multicolumn{1}{l}{1.52E-01}
      \\
    \bottomrule
    \end{tabular}%
  \label{tab1}%
\end{table}%
Table \ref{tab1} and Figure \ref{fig1} show the difference of solutions among  (ori), (red), and (ie) when the number of particles is $10^4$ and the radius of each particle is $10^{-4}$ cm with various values for the distance $d$.  For $10^4$ particles, the error $e$ is smallest, equal to 1.29\%, when $d$ is $5\times 10^{-2}$ cm, or $a/d=2\times 10^{-3}$. The error grows slowly when $d$ is slightly away from this point, and it is greater than 5\% when $d \ge 5.8\times 10^{-2}$ cm or $d \le 4.2\times 10^{-2}$ cm. The error is less than 10\% when $1.6\times 10^{-3} \le \frac{a}{d} \le 2.5\times 10^{-3}$. The solutions to the reduce system and the integral equation are very closed since the reduce system is essentially the Riemann sum of the integral equation.

\begin{table}[H]
  \centering
  \tabcolsep=0.11cm
  \small
  \caption{Solution comparison of (ori), (red), and (ie) with $M=10^5$, $a=10^{-5}$, and different $d$.}
    \begin{tabular}{crrrrrr}
    \toprule
    \multicolumn{7}{c}{M=1.00E+5, a=1.00E-5}
      \\
    \midrule
    d &
      \multicolumn{1}{l}{1.00E-02} &
      \multicolumn{1}{l}{1.50E-02} &
      \multicolumn{1}{l}{2.00E-02} &
      \multicolumn{1}{l}{2.30E-02} &
      \multicolumn{1}{l}{2.50E-02} &
      \multicolumn{1}{l}{3.00E-02}
      \\
    a/d &
      \multicolumn{1}{l}{1.00E-03} &
      \multicolumn{1}{l}{6.67E-04} &
      \multicolumn{1}{l}{5.00E-04} &
      \multicolumn{1}{l}{4.35E-04} &
      \multicolumn{1}{l}{4.00E-04} &
      \multicolumn{1}{l}{3.33E-04}
      \\
    (ori) vs. (red) &
      \multicolumn{1}{l}{9.04E-02} &
      \multicolumn{1}{l}{5.34E-02} &
      \multicolumn{1}{l}{1.64E-02} &
      \multicolumn{1}{l}{1.21E-02} &
      \multicolumn{1}{l}{2.74E-02} &
      \multicolumn{1}{l}{6.57E-02}
      \\
    (ie) vs. (ori) &
      \multicolumn{1}{l}{1.02E-01} &
      \multicolumn{1}{l}{6.89E-02} &
      \multicolumn{1}{l}{3.56E-02} &
      \multicolumn{1}{l}{1.56E-02} &
      \multicolumn{1}{l}{1.44E-02} &
      \multicolumn{1}{l}{5.29E-02}
      \\
    (ie) vs. (red) &
      \multicolumn{1}{l}{3.04E-03} &
      \multicolumn{1}{l}{3.04E-03} &
      \multicolumn{1}{l}{3.04E-03} &
      \multicolumn{1}{l}{3.04E-03} &
      \multicolumn{1}{l}{3.04E-03} &
      \multicolumn{1}{l}{3.04E-03}
      \\
    Error sum $e$ &
      \multicolumn{1}{l}{1.96E-01} &
      \multicolumn{1}{l}{1.25E-01} &
      \multicolumn{1}{l}{5.51E-02} &
      \multicolumn{1}{l}{3.07E-02} &
      \multicolumn{1}{l}{4.49E-02} &
      \multicolumn{1}{l}{1.22E-01}
      \\
    \bottomrule
    \end{tabular}%
  \label{tab2}%
\end{table}%

Table \ref{tab2} and Figure \ref{fig2} show the difference of solutions among  (ori), (red), and (ie) when the number of particles is $10^5$, and the radius of a particle is $10^{-5}$ cm with different values for the distance $d$. In this case, the error $e$ is smallest, equal to 3\%, when $d$ is $2.3\times 10^{-2}$ cm, or $a/d=4.35\times 10^{-4}$. The error grows quite slowly when $d$ increases or decreases from this point. The error is less than 10\% when $3.7 \times 10^{-4} \le \frac{a}{d} \le 6\times 10^{-4}$.

\begin{table}[H]
  \centering
  \tabcolsep=0.11cm
  \small
  \caption{Solution comparison of (ori), (red), and (ie) with $M=10^6$, $a=10^{-6}$, and different $d$.}
    \begin{tabular}{crrrrrr}
    \toprule
    \multicolumn{7}{c}{M=1.00E+6, a=1.00E-6}
      \\
    \midrule
    d &
      \multicolumn{1}{l}{5.00E-03} &
      \multicolumn{1}{l}{7.00E-03} &
      \multicolumn{1}{l}{9.00E-03} &
      \multicolumn{1}{l}{9.50E-03} &
      \multicolumn{1}{l}{1.00E-02} &
      \multicolumn{1}{l}{1.50E-02}
      \\
    a/d &
      \multicolumn{1}{l}{2.00E-04} &
      \multicolumn{1}{l}{1.43E-04} &
      \multicolumn{1}{l}{1.11E-04} &
      \multicolumn{1}{l}{1.05E-04} &
      \multicolumn{1}{l}{1.00E-04} &
      \multicolumn{1}{l}{6.67E-05}
      \\
    (ori) vs. (red) &
      \multicolumn{1}{l}{8.26E-02} &
      \multicolumn{1}{l}{5.00E-02} &
      \multicolumn{1}{l}{1.73E-02} &
      \multicolumn{1}{l}{9.09E-03} &
      \multicolumn{1}{l}{1.62E-03} &
      \multicolumn{1}{l}{8.08E-02}
      \\
    (ie) vs. (ori) &
      \multicolumn{1}{l}{8.95E-02} &
      \multicolumn{1}{l}{5.41E-02} &
      \multicolumn{1}{l}{1.86E-02} &
      \multicolumn{1}{l}{9.77E-03} &
      \multicolumn{1}{l}{9.16E-04} &
      \multicolumn{1}{l}{8.76E-02}
      \\
    (ie) vs. (red) &
      \multicolumn{1}{l}{3.04E-03} &
      \multicolumn{1}{l}{3.04E-03} &
      \multicolumn{1}{l}{3.04E-03} &
      \multicolumn{1}{l}{3.04E-03} &
      \multicolumn{1}{l}{3.04E-03} &
      \multicolumn{1}{l}{3.04E-03}
      \\
    Error sum $e$ &
      \multicolumn{1}{l}{1.75E-01} &
      \multicolumn{1}{l}{1.07E-01} &
      \multicolumn{1}{l}{3.89E-02} &
      \multicolumn{1}{l}{2.19E-02} &
      \multicolumn{1}{l}{5.58E-03} &
      \multicolumn{1}{l}{1.71E-01}
      \\
    \bottomrule
    \end{tabular}%
  \label{tab3}%
\end{table}%

Table \ref{tab3} and Figure \ref{fig3} show the difference of solutions among  (ori), (red), and (ie) when the number of particles is $10^6$, the radius of a particle is $10^{-6}$ cm, and the distance $d$ varies. In this case, the error $e$ is smallest, equal to 0.56\%, when $d$ is $1\times 10^{-2}$ cm, that is $a/d=1\times 10^{-4}$. The error grows slightly when $d$ is between $7\times 10^{-3}$ cm and $1.5\times 10^{-2}$ cm. After that, the error increases significantly. The error is less than 10\% when $ 8\times 10^{-5} \le \frac{a}{d} \le 1.4\times 10^{-4}$.

Next, we will look at the best ratio $a/d$ for each radius $a$ for which the error sum $e$ is smallest, i.e. the best ratio $a/d$ will minimize the solution differences among the (ori), (red) and (ie). The error sum is used as the criterion for the optimization. For each number of particles $M$ and radius $a$, we feed many different values of $d$ to find the smallest error sum.

Note that since we use uniform distribution and the size of the domain $\Omega$, radius $a$ and number of particles $M$ are fixed, size of $\Omega$ is 1 cm, $a$ and $M$ are fixed in each test case, we cannot increase the distance between neighboring particles $d$ to the size of the cube $\Omega$ or decrease $d$ to be less than $a$. The distance $d$ must be of order $O\left(a^{\frac{2-\kappa}{3}}\right)$ as described in \eqref{eq11} so that all the particles lie in the domain $\Omega$.

\begin{table}[htbp]
  \centering
  \tabcolsep=0.11cm
  \small
  \caption{The best ratio $a/d$ for each radius $a$.}
    \begin{tabular}{clll}
    \toprule
    M &
      1.00E+06 &
      1.00E+05 &
      1.00E+04
      \\
    \midrule
    a &
      1.00E-06 &
      1.00E-05 &
      1.00E-04
      \\
    d &
      1.00E-02 &
      2.30E-02 &
      5.00E-02
      \\
    a/d &
      1.00E-04 &
      4.35E-04 &
      2.00E-03
      \\
    (ori) vs. (red) &
      1.62E-03 &
      1.21E-02 &
      6.52E-03
      \\
    (ie) vs. (ori) &
      9.16E-04 &
      1.56E-02 &
      4.57E-03
      \\
    (ie) vs. (red) &
      3.04E-03 &
      3.04E-03 &
      1.83E-03
      \\
    Error sum $e$ &
      5.58E-03 &
      3.07E-02 &
      1.29E-02
      \\
    \bottomrule
    \end{tabular}%
  \label{tab4}%
\end{table}%
\begin{figure}[htbp]
    \centering
    \includegraphics[scale=1.1]{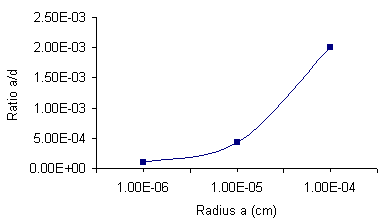}
    \caption{The best ratio $a/d$ for each radius $a$.}
    \label{fig4}
\end{figure}

Table \ref{tab4} and Figure \ref{fig4} show the best ratios $\frac{a}{d}$, corresponding to the smallest error sums, when $a$ is $10^{-4}$, $10^{-5}$ and $10^{-6}$ cm, and $M$ is $10^{4}, 10^{5}$ and $10^{6}$ particles, respectively. For instance, the best ratio $a/d$ at $a=10^{-6}$ cm is $1\times 10^{-4}$. The optimal values of $d$ for which the optimizations, the smallest error sums, are obtained are also given. As one can see, the optimal value of $d$ is within a small finite range and depends on the radius $a$. As $a$ gets smaller, this range becomes smaller as well. The quality of the approximation of the solution to the wave scattering problem depends on this range.

\section{Conclusions} \label{sec5}
The numerical experiment shows that the errors, i.e. solution differences of  (ori), (red), and (ie), depend greatly on the radius of particles, $a$, the number of particles, $M$, and the distance between neighboring particles, $d$. The numerical results help us to better understand the asymptotic solutions to the problem of acoustic wave scattering by many small impedance particles and the possibility of creating materials with any desired refraction coefficient by using the asymptotic approach. Indeed, for acoustic wave scattering, there is an optimal value of the ratio $\frac{a}{d}$ for which the error is acceptable and the asymptotic solution to  (red) can be used as a good approximation to the solutions of  (ori) and (ie). This would help to simplify the computation process immensely, specifically when the number of particles is extremely large and the radius of particles is very small.

In the future, we will consider developing a new algorithm for conducting the experiment with a larger number of particles, say $M$ from $10^7$ up to $10^{12}$. The current algorithm does not allow us to go beyond $10^6$ particles since it requires $O(n^2)$ operations for matrix-vector multiplication in the iterative process, which is very expensive in terms of computation time.  \\

\noindent {\bf Acknowledgements.} The author is grateful to Professor A. G. Ramm for teaching wave scattering and Professor N. Albin for sharing knowledge about high performance computing for doing the experiments.

The computing for this project was performed on the Beocat Research Cluster at Kansas State University.

\bibliographystyle{ieeetr}

\end{document}